\documentclass[11pt]{amsart}

\usepackage{fullpage}

\usepackage{url}

\usepackage{amssymb}

\usepackage{cite}

\renewcommand{\a}{\alpha}

\newcommand{\g}{\gamma}
\renewcommand{\d}{\delta}

\newcommand{\e}{\varepsilon}

\renewcommand{\S}{\Sigma}

\renewcommand{\l}{\lambda}

\newcommand{\cO}{{\mathcal O}}

\newcommand{\cC}{{\mathcal C}}
\newcommand{\cM}{{\mathcal M}}

\newcommand{\cE}{{\mathcal E}}

\newcommand{\cV}{{\mathcal V}}

\newcommand{\cN}{{\mathcal N}}

\newcommand{\cH}{{\mathcal H}}

\newcommand{\cI}{{\mathcal I}}

\newcommand{\bR}{\mathbb R}

\newcommand{\bZ}{\mathbb Z}

\newcommand{\bE}{\mathbb E}

\newcommand{\be}{\begin{equation}}
\newcommand{\ee}{\end{equation}}

\newcommand{\bel}[1]{\begin{equation}\label{#1}}
\newcommand{\beaa}{\begin{eqnarray*}}
\newcommand{\bea}{\begin{eqnarray}}
\newcommand{\beal}[1]{\begin{eqnarray}\label{#1}}
\newcommand{\bean}{\begin{eqnarray}\nonumber}
\newcommand{\beadl}[1]{\begin{deqarr}\label{#1}}
\newcommand{\eeadl}[1]{\arrlabel{#1}\end{deqarr}}
\newcommand{\eeal}[1]{\label{#1}\end{eqnarray}}
\newcommand{\eead}[1]{\end{deqarr}}
\newcommand{\eea}{\end{eqnarray}}
\newcommand{\eeaa}{\end{eqnarray*}}

\renewcommand{\to}{\rightarrow}

\renewcommand{\phi}{\varphi}
\renewcommand{\epsilon}{\varepsilon}
\renewcommand{\hat}{\widehat}

\newcommand{\dm}{{\partial M}}

\newcommand{\w}{\widetilde}

\theoremstyle{plain}
\newtheorem{theorem}{Theorem}[section]

\newtheorem{lemma}[theorem]{Lemma}

\newtheorem{proposition}[theorem]{Proposition}

\theoremstyle{definition}

\def\endproof{\qed \medskip}
\def\blacksquare{\hbox to .60em {\vrule width .60em height .60em}}

\numberwithin{equation}{section}

\date{\today}

\begin{document}

\title[]{On the Bartnik conjecture for the \\ static vacuum Einstein equations}

\author[]{Michael T. Anderson}

\address{Dept.~of Mathematics, Stony Brook University, Stony Brook, NY 11790}
\email{anderson@math.sunysb.edu}

\thanks{The author is partially supported by NSF Grant DMS 1205947  \\
PACS 2010: 04.20.-q, 04.20.Cv, 02.40.Vh, 02.30.Jr}

\abstract{We prove that given any smooth metric $\g$ and smooth positive function $H$ on $S^{2}$, there 
is a constant $\l > 0$, depending on $(\g, H)$, and an asymptotically flat solution $(M, g, u)$ of the 
static vacuum Einstein equations on $M = \bR^{3}\setminus B^{3}$, such that the induced metric 
and mean curvature of $(M, g, u)$ at $\dm$ are given by $(\g, \l H)$. This gives a partial resolution of 
a conjecture of Bartnik.}
\endabstract

\maketitle

\setcounter{section}{0}
\setcounter{equation}{0}

\section{Introduction}

  The Bartnik conjecture \cite{B2}, \cite{B4} is a statement concerning the existence and uniqueness of solutions 
to a geometrically and physically natural boundary value problem for the static Einstein equations in vacuum. 
In general relativity, it is important in connection with the concepts of mass and quasi-local mass of 
initial data sets (Cauchy surfaces). Mathematically, the conjecture concerns the unique global solvability 
of a very interesting but rather complicated elliptic boundary value problem for a system of PDE's 
on a $3$-manifold $M$. 

\medskip 

  We first recall the statement of the conjecture. Let $B^{3} \subset \bR^{3}$ be a 3-ball (not necessarily round) 
and let $M$ be a 3-manifold (with boundary) diffeomorphic to $\bR^{3}\setminus B^{3}$. Thus $\dm$ is 
diffeomorphic to the 2-sphere $S^{2}$. The free data for a static vacuum solution on $M$ are a pair 
$(g, u)$ consisting of a smooth asymptotically flat Riemannian metric $g$ on $M$ and a smooth positive 
(potential) function $u: M \to \bR^{+}$ asymptotic to $1$ at infinity. The static vacuum Einstein equations 
on $(g, u)$ are given by 
\begin{equation}\label{stat}
uRic_{g} = D^{2}u, \ \ \Delta u = 0,
\end{equation}
where the Hessian $D^{2}$ and Laplacian $\Delta = tr D^{2}$ are taken with respect
to $g$. From the metric $g$ and potential $u$, one can construct a 4-dimensional solution to 
the (usual) vacuum Einstein equations 
\begin{equation}\label{vac}
Ric_{g_{\cM}} = 0,
\end{equation}
where $\cM = \bR \times M$ and the 4-metric is given by $g_{\cM} = \pm u^{2}dt^{2} + g$. 
The equations \eqref{stat} on $M$ and \eqref{vac} on $\cM$ are equivalent. 

\medskip 

  Given a solution $(M, g, u)$ of the static vacuum Einstein equations, let $\g$ be the 
Riemannian metric induced on $S^{2} \simeq \dm$ and let $H$ be the mean curvature of 
$\dm \subset (M, g)$, (with respect to the inward unit normal into $M$). View the pair $(\g, H)$ 
as boundary data of the static vacuum solution $(M, g, u)$ and suppose the data $(g, u)$ are 
$C^{\infty}$ smooth up to $\dm$. The boundary data thus takes values in the space 
$Met^{\infty}(S^{2})\times C^{\infty}(S^{2})$ of smooth metrics and smooth functions on 
$S^{2}$. Observe that the potential $u$ on $\dm$ is not part of the boundary data. Let 
$C_{+}^{\infty}(S^{2})$ be the space of positive $C^{\infty}$ functions $H: S^{2} \to \bR^{+}$. 

  The Bartnik conjecture then states that arbitrary boundary data $(\g, H) \in Met^{\infty}(S^{2})
\times C_{+}^{\infty}(S^{2})$ are realized uniquely (up to diffeomorphism) by an asymptotically 
flat solution $(M, g, u)$ of the static vacuum Einstein equations \eqref{stat} on $M$, smooth up to 
$\dm$. Thus, given an arbitrary smooth metric $\g$ on $S^{2}$ and an arbitrary positive smooth function 
$H$ on $S^{2}$, there exists a smooth asymptotically flat solution $(M, g, u)$ of the static vacuum Einstein 
equations for which the induced metric and mean curvature at $\dm$ are given by $(\g, H)$. 
Further, the solution is unique, up to a smooth diffeomorphism fixing $\dm$. 

\smallskip 

  In this paper, we obtain a partial resolution of the existence part of the Bartnik conjecture. One  
version of the main result is the following. 
\begin{theorem}
Given any smooth boundary data $(\g, H) \in Met^{\infty}(S^{2})\times C_{+}^{\infty}(S^{2})$, 
there is a constant $\l > 0$, depending on $(\g, H)$, and an asymptotically flat static vacuum 
solution $(M, g, u)$, smooth up to $\dm$, such that the induced metric and mean curvature 
of $(M, g, u)$ at $\dm$ are given by $(\g, \l H)$. 
\end{theorem}

   Thus, up to a single (multiplicative) scalar degree of indeterminacy on the mean curvature, 
arbitrary smooth data $(\g, H)$ are realized as the boundary metric and mean curvature of a smooth 
asymptotically flat static vacuum solution $(M, g, u)$. In particular, the boundary metric $\g$ on 
$S^{2}$ may be arbitrarily specified. Clearly (by rescaling $g$) one may instead shift the indeterminacy 
to a scaling of the boundary metric, (so that given $(\g, H)$, there is a constant $\l > 0$ such that 
$(\l \g, H)$ are realized as boundary data of an asymptotically flat static vacuum solution). 

\medskip 

  The idea of considering rescalings $\l H$ of the mean curvature appears in earlier work of 
Jauregui \cite{J}, (cf.~also \cite{JMT}), in connection with filling {\it in} prescribed boundary data $(\g, \l H)$ on 
$S^{2}$ with a metric of non-negative scalar curvature on a $3$-ball $B^{3}$ with boundary $S^{2}$. 
Similarly, it appears in the context of filling {\it in} such boundary data by static vacuum solutions 
(the ``interior Bartnik problem") in \cite{An2}. 

   In the interior context, such a rescaling of $H$ is 
{\it necessary}; for a typical $\bR^{+}$-orbit $(\g, \l H)$, $\l \in \bR^{+}$, only for at most 
finitely many $\l$ are the data $(\g, \l H)$ realized as boundary data of interior static vacuum solutions. 
In fact, the space of interior (strictly) static vacuum solutions is generically transverse to the $\bR^{+}$-orbits 
of $H$, cf.~\cite{An2}. 

  Such a restriction is not necessary in the current ``exterior" context. There are simple examples 
where for a given data $(\g, H)$, all $H$-rescalings $(\g, \l H)$ are also realized as boundary data of 
static vacuum solutions. The simplest instance of this is the curve of Schwarzschild metrics, given as 
\be \label{Sch}
g_{m} = (1 - \frac{2m}{r})^{-1}dr^{2} + r^{2}g_{S^{2}(1)},
\ee
on the space-like hypersurface $M \simeq \bR^{3}\setminus B$ where $r \geq 1 > 0$ and $r > 2m$. 
Assume $m < \frac{1}{2}$. Then the boundary $\dm$ may be taken at the locus $r = 1$ and the metric on 
$\{r = 1\}$ is fixed, 
$$\g = \g_{S^{2}(1)}.$$
A simple computation shows that at the boundary $\dm = \{r = 1\}$,
$$H = 2\sqrt{1 - 2m},$$
so that $H$ is a constant depending only on $m$. As $m$ varies over the interval $(-\infty, \frac{1}{2})$, 
$H \in (0, \infty)$ is monotone decreasing. 

\medskip 

  It is worthwhile to place the Bartnik conjecture and Theorem 1.1 in a more geometric setting. 
To do this, let $\bE^{m,\a}$ be the space of all asymptotically flat (AF) solutions $(g, u)$ of the 
static vacuum equations \eqref{stat}, which are $C^{m,\a+}$ on $M$ up to $\dm$. Here $C^{m,\a+}$ 
is the usual ``little" H\"older space of functions, i.e.~the closure of the space of $C^{\infty}$ smooth 
functions on $M$ with respect to the $C^{m,\a}$ H\"older norm. It is well-known, cf.~\cite{O} for instance, 
that $C^{m,\a+}$ is a closed separable subspace of the H\"older space $C^{m,\a}$ (which itself is not separable). 
We assume throughout the paper that $m \geq 3$ and $\a \in (0,1)$. As above, we also assume $u \to 1$ 
at infinity in $M$. The topology on $\bE^{m,\a}$ is the standard weighted $C^{m,\a}$ H\"older topology used 
for AF manifolds, cf.~\cite{AK} for instance. 

   Let $\mathrm{Diff}_{1}^{m+1,\a}(M)$ be the group of $C^{m+1,\a+}$ diffeomorphisms 
$\psi: M \to M$ with $\psi = id$ on $\dm$. The group $\mathrm{Diff}_{1}^{m+1,\a}(M)$ acts freely on 
$\bE^{m,\a}$ fixing the boundary data $(\g, H)$, (since an isometry of $(M, g)$ fixing $\dm$ pointwise 
is necessarily the identity). The quotient space $\cE^{m,\a}$ is the moduli space of static vacuum solutions; 
an element in $\cE^{m,\a}$ represents an equivalence class of isometric metrics. 

 Next, as above abusing notation slightly, let $Met^{m,\a}(S^{2})$ be the space of $C^{m,\a+}$ metrics on 
$S^{2} \simeq \dm$ and $C^{m-1,\a}(S^{2})$ be the space of $C^{m-1,\a+}$ functions on $S^{2}$. 
One thus has a natural map, mapping a static vacuum solution to its Bartnik boundary data:
\begin{equation}\label{PiB}
\Pi_{B}:  \cE^{m,\a} \to Met^{m,\alpha}(S^{2})\times C^{m-1,\a}(S^{2}),
\end{equation}
$$\Pi_{B}(g) = (\g, H).$$

It is proved in \cite{An1}, \cite{AK} that the space $\cE^{m,\a}$ is a smooth (infinite dimensional) Banach 
manifold, and the map $\Pi_{B}$ is $C^{\infty}$ smooth and Fredholm, of Fredholm index 0.\footnote{This 
is proved in \cite{An1}, \cite{AK} with respect to the usual H\"older $C^{m,\a}$ spaces. The proof relies only 
on the fact that the linearized static vacuum equations, together with the linearization of the boundary conditions $(\g, H)$, 
form a formally self-adjoint elliptic boundary value problem in a suitable gauge, (cf.~Section 3 and in particular 
Theorem 3.5 of [3]). The same proof, based on elliptic regularity and Fredholm theory, applies also to the little 
H\"older spaces $C^{m,\a+}$ and the Banach manifolds modeled thereon. While such extensions of elliptic 
regularity and Fredholm theory to the little H\"older spaces are forklore, an explicit and general proof is given in Appendix 
A of \cite{D} for instance.} This means that the derivative or linearization $D\Pi_{B}$, at any solution $(M, g, u) \in \cE^{m,\a}$, 
has finite dimensional kernel $Ker D\Pi_{B}$ and finite dimensional cokernel $Coker D\Pi_{B}$, 
$$dim Ker D\Pi_{B} < \infty \ \ {\rm and} \ \  dim Coker D\Pi_{B} < \infty,$$
and the range of $D\Pi_{B}$ is closed. Further the Fredholm index, given by the difference of these 
terms, satisfies 
\be \label{ind0}
dim Ker D\Pi_{B} - dim Coker D\Pi_{B} = 0.
\ee
This result implies that the boundary map $\Pi_{B}$ is at least locally reasonably well-behaved, as a 
map of Banach manifolds. 

  Let $\cE_{+}^{m,\a}$ be the open submanifold of $\cE^{m,\a}$ for which
the mean curvature $H$ on $\dm$ is positive, i.e.
$$\cE_{+}^{m,\a} = (\Pi_{B})^{-1}(Met^{m,\a}(S^{2})\times C_{+}^{m-1,\a}(S^{2})).$$
The Bartnik conjecture above may thus be rephrased as stating that the smooth map $\Pi_{B}$,
restricted to the open submanifold $\cE_{+}^{m,\a}$,
\begin{equation}\label{PiB+}
\Pi_{B}: \cE_{+}^{m,\a} \to Met^{m,\a}(S^{2})\times C_{+}^{m-1,\a}(S^{2}),
\end{equation}
is surjective and injective. (Here we generalize slightly from $m = \infty$ to general $m \geq 3$). 
If one adds the condition that the linearization $D\Pi_{B}$ is an isomorphism everywhere (i.e.~$Ker 
D\Pi_{B} = 0$ everywhere) it follows from the inverse function theorem in Banach manifolds that 
$\Pi_{B}$ is conjecturally a smooth global diffeomorphism. 

   However, it is proved in \cite{AK} that this most optimistic version of the conjecture does not hold; 
$\Pi_{B}$ in \eqref{PiB+} cannot be a diffeomorphism. For example, consider the ``trivial" 
static solution $\bR^{3}$ with $u = 1$ and $g = g_{Eucl}$. It is easy to construct a curve of embedded 
positive mean curvature $2$-spheres $F_{t}: S^{2} \to \bR^{3}$, $t \in (0,1]$ such that as $t \to 0$, 
the boundary data converge to a limit in the target space, i.e.~$(\g_{t}, H_{t}) \to (\g, H) 
\in Met^{m,\a}(S^{2})\times C_{+}^{m-1,\a}(S^{2})$ with $F_{t} \to F_{0}$, where $F_{0}(S^{2})$ is 
only an immersed but not an embedded 2-sphere. It follows that the exterior or unbounded 
components of $\bR^{3} \setminus F_{t}(S^{2})$ define static vacuum solutions $(M_{t}, 
(g_{Eucl})_{t}, 1)$ for which the boundary data $\Pi_{B}(M_{t}, (g_{Eucl})_{t}, 1)$ converge to a limit in 
$Met^{m,\a}(S^{2})\times C_{+}^{m-1,\a}(S^{2})$, but although the limit $(M_{0}, (g_{Eucl})_{0}, 1)$ makes 
sense, it is not a static vacuum solution on a smooth manifold with boundary. The structure of a 
manifold-with-boundary degenerates, with uniform control on the boundary data $(\g, H)$. This 
means that the inverse map $(\Pi_{B})^{-1}$ (assuming it exists) cannot be continuous everywhere. 

   The same construction can be carried out in any background (exterior) static vacuum solution $(M, g, u)$ 
(in place of $(\bR^{3}, g_{Eucl}, 1))$ so that this phenomenon is ubiquitous.  

\medskip 

  To understand Theorem 1.1 in this context, define an equivalence relation on $C_{+}^{m-1,\a}(S^{2})$ 
by $[H_{1}] = [H_{2}]$ if and only if 
\be \label{equiv}
H_{2} = \l H_{1},
\ee
for some constant $\l > 0$. Let $\cC_{+}^{m-1,\a}(S^{2}) = C_{+}^{m-1,\a}(S^{2})/ \sim$ 
be the space of equivalence classes; this is the quotient of $C_{+}^{m-1,\a}(S^{2})$ by the 
$\bR^{+}$-action given by multiplication by a positive scalar. Since the action is smooth and free, 
the space $\cC_{+}^{m-1,\a}(S^{2})$ is also a smooth Banach manifold, and so one has an induced 
smooth map 
\be \label{Pi}
\Pi: \cE_{+}^{m,\a} \to Met^{m,\a}(S^{2}) \times \cC_{+}^{m-1,\a}(S^{2}),
\ee
$$\Pi(M, g, u) = (\g, [H]).$$
The map $\Pi$ is a smooth Fredholm map, of Fredholm index 1, i.e.
\be \label{ind1}
dim Ker D\Pi - dim Coker D\Pi = 1.
\ee
Theorem 1.1 then follows from the following slightly more general statement. 

\begin{theorem}
The map $\Pi$ in \eqref{Pi} is surjective. 
\end{theorem}

   The proof of Theorems 1.1 and 1.2 follow in Section 2 below. The proof builds on the prior work of the author in 
\cite{AK} which will be used quite often here. In particular, we refer to \cite{AK} for a more detailed presentation of 
certain definitions. Section 3 contains a discussion of related work and lists several open problems suggested 
by Theorem 1.2. 

  I am most grateful to the referee and to Piotr Chru\'sciel for very useful comments leading to a substantial 
improvement of the original version of this paper.

\section{Proofs.}
\setcounter{equation}{0}

  Given the results of \cite{AK}, the primary difficulty in understanding 
the validity of the Bartnik conjecture is the degeneration of the manifold-with-boundary structure under 
control of the boundary data $(\g, H)$ discussed in the Introduction. 

   Given a static vacuum solution $(M, g, u) \in \cE^{m,\a}$, let $N$ be the unit inward normal of $\dm$ in $M$ 
and let $inj_{\dm}$ denote the injectivity radius of the normal exponential map from $\dm$ into $(M, g, u)$. 
The normal exponential map $exp_{N}: N(\dm) \to M$ maps vectors in the normal bundle of $\dm$ into $M$ 
and $inj_{\dm}$ is the largest value such that $exp_{N}$ is a diffeomorphism when restricted to normal vectors 
$aN$ with $a < inj_{\dm}$.

   The behavior described in the paragraph following \eqref{PiB+} is characterized by the fact that there are curves 
(or sequences) of static vacuum solutions $(M, g_{t}, u_{t})$ (with uniformly bounded curvature and derivatives of curvature for 
instance) such that 
$$inj_{\dm_{t}} \to 0, \ \ {\rm as} \ \ t \to 0,$$
but with $(\g_{t}, H_{t}) \to (\g, H) \in Met^{m,\a}(S^{2})\times C_{+}^{m-1,\a}(S^{2})$. Equivalently, there 
are a pair of curves $x_{t}^{1}$, $x_{t}^{2}$ in $\dm_{t}$ with $dist_{\dm_{t}}(x_{t}^{1}, x_{t}^{2}) \geq c_{0} > 0$
but such that 
$$dist_{(M,g_{t})}(x_{t}^{1}, x_{t}^{2}) \to 0, \ \ {\rm as} \ \ t \to 0.$$
Thus the distance between $x_{t}^{1}$ and $x_{t}^{2}$ remains bounded away from zero in $\dm$, but converges 
to zero in $(M, g_{t})$. 

\medskip 

  Motivated by the work in \cite{AK}, we will use a somewhat stronger notion than the injectivity radius. 
Let $S(t) = \{x\in M: dist(x, \dm) = t\}$ be the geodesic $t$-sphere about $\dm$ and 
$B(t) = \{x \in M: dist(x, \dm) < t\}$ the associated open geodesic ball about $\dm$. 
For $t < inj_{\dm}$, $S(t)$ is a $C^{m+1,\a}$ smooth hypersurface, diffeomorphic to $\dm$, with 
$\partial B(t) = S(t) \cup \dm$. (For $t > inj_{\dm}$, $S(t)$ is no longer in general smooth). 

  Define the {\em outer-minimizing radius} of $\dm$ in $(M, g, u) \in \cE_{+}^{m,\a}$ by 
\be \label{h}
\cO_{\dm} = \sup \{\rho: area(\S) > area(\dm)\},
\ee
where $\S$ is any smooth surface in $B(\rho)$ homologous to $\dm$, $\S \neq \dm$. Thus $\cO_{\dm}$ 
is the largest radius such that $\dm$ is outer-minimizing in $B(\cO_{\dm})$, i.e.~$\dm$ has least area 
among all surfaces $\S \subset B(\cO_{\dm})$ homologous to $\dm$. Recall that the concept of 
outer-minimizing boundary plays an important role in \cite{AK}.

  By definition, $H > 0$ on $\dm$ for $(M, g, u) \in \cE_{+}^{m,\a}$. Hence, there is an 
$\e > 0$ (depending on $(M, g, u)$) such that $H(s) > 0$ for all $s < \e$, where $H(s)$ is the 
mean curvature of $S(s)$. Thus, $B(\e)$ has a foliation by surfaces of positive mean curvature. 
Since $H > 0$ corresponds to infinitesimal area expansion, it is easy to see (and well-known) that 
$\dm$ is outer-minimizing in $B(\e)$, so that 
$$\cO_{\dm} \geq \e > 0.$$ 
Thus $\cO_{\dm}$ is a positive function on $\cE_{+}^{m,\a}$.

  It is also well-known that if the Riemann curvature $Rm$ of $(M, g)$ is uniformly bounded, $|Rm| \leq \Lambda$, 
then $\cO_{\dm} \geq d_{0} > 0$ implies $inj_{\dm} \geq d_{1} > 0$, where $d_{1}$ depends only on $d_{0}$ and 
$\Lambda$; (for a proof see the beginning of the proof of [3, Theorem 4.1] for instance). The converse is also true, 
i.e.~ $inj_{\dm} \geq d_{1} > 0$ implies $\cO_{\dm} \geq d_{0} > 0$, with $d_{0} = d_{0}(d_{1}, \Lambda)$. 

\medskip 

  The following result from \cite{AK} shows that there is no degeneration in $\cE_{+}^{m,\a}$ with uniform 
control on the boundary data $(\g, H)$ and a uniform lower bound on $\cO_{\dm}$. 

\begin{proposition}
Let $S$ be a set of static vacuum metrics $(M, g, u) \in \cE_{+}^{m,\a}$ for which the corresponding set of boundary 
data $(\g, H)$ is a compact subset of $Met^{m,\a}(S^{2})\times C_{+}^{m-1,\a}(S^{2})$ and for which there exists 
a constant $d_{0} > 0$ such that for all $(M, g, u) \in S$,  
$$\cO_{\dm} \geq d_{0} > 0.$$
Then $S$ is compact in $\cE_{+}^{m,\a}$. 
\end{proposition}

{\bf Proof:} This result is discussed following the statement of [3, Theorem 1.2] and the proof, which follows from 
the proof of [3, Theorem 4.3], is discussed in [3, Remark 4.4]. 

{\endproof}

  The basic idea in the proof of Theorem 1.2 is to build into the boundary data $(\g, H)$ a term 
controlling the behavior of $\cO_{\dm}$. We do this by modifying the boundary map $\Pi_{B}$ by a suitable term 
related to $\cO_{\dm}$. Of course this is no longer a ``boundary map" in the strict sense, but as will be 
seen, it is sufficient to prove Theorem 1.2. 
   
   We first note the following: 
   
\begin{lemma}  
The function
$$\cO_{\dm}: \cE_{+}^{m,\a} \to \bR^{+} \cup \{\infty\}$$
is lower semi-continuous. 
\end{lemma}

{\bf  Proof:}   Let $(M_{i}, g_{i}, u_{i})$ be a sequence in $\cE_{+}^{m,\a}$ converging to a limit $(M, g, u)$. 
Assume that
$$\cO_{0} \equiv \liminf \cO_{\dm}(g_{i}) < \infty,$$
(since otherwise there is nothing to prove). In each closed ball $\bar B_{i}(\cO_{\dm}(g_{i}))$ about $\dm$ in $(M, g_{i})$ 
there is a surface $\S_{i} \neq \dm$ such that 
$$area_{g_{i}}(\S_{i}) = area_{g_{i}}(\dm).$$
Since the convergence to the limit $(M, g, u)$ is sufficiently smooth (i.e.~$C^{m,\a}$), there exists a limit surface 
$\S \subset \bar B(\cO_{0}) \subset (M, g)$ such that 
$$area_{g}(\S) = area_{g}(\dm).$$
It follows that 
$$\lim \inf \cO_{\dm}(g_{i}) \geq \cO_{\dm},$$
which proves lower semi-continuity. 

{\endproof}

  Let $\cV^{m,\a} \subset \cE_{+}^{m,\a}$ be the set of static vacuum solutions $(M, g, u)$ such that 
$\dm$ is outer-minimizing in all of $M$, i.e.
$$\cV^{m,\a} = \{(M, g, u) \in \cE_{+}^{m, \a}: \cO_{\dm} = \infty\}.$$
We claim that $\cV^{m,\a}$ is open in $\cE_{+}^{m,\a}$. This follows from the fact that the solutions $(M, g, u)$ 
are asymptotically flat, (and the convergence in $\cE^{m,\a}$ preserves this structure). Namely, surfaces 
$\S \subset M \setminus B(R)$ homologous to $\dm$ have very large area (for $R$ sufficiently large) so 
that for all such surfaces 
$$area(\S) > area(\dm).$$ 
In particular if $(M, g_{i}, u_{i}) \to (M, g, u) \in \cE^{m,\a}$ and $\cO_{\dm}(M) = \infty$, then $\cO_{\dm}(g_{i}) = 
\infty$ for $i$ sufficiently large. 

  The reciprocal function 
\be \label{reci}
\cI_{\dm} = \frac{1}{\cO_{\dm}}: \cE_{+}^{m,\a} \to \bR^{+} \cup \{0\} \subset \bR
\ee
is upper semi-continuous. (Recall $f$ is upper semi-continuous if $\limsup_{x \to x_{0}} f(x) \leq f(x_{0})$). 

    Now the space $\cE^{m,\a}$ is a smooth separable Banach manifold locally modeled on the little 
H\"older spaces $C^{m,\a+}$ and has smooth partitions of unity, cf.~\cite{FM}, \cite{BFT} 
and references therein. In particular $\cE^{m,\a}$ is a perfectly normal topological space and 
by a classical theorem of Baire, any upper semi-continuous function on $\cE_{+}^{m,\a}$ is then the 
pointwise limit of a decreasing sequence of continuous functions; see for instance \cite{Mu} for the 
definition and discussion of such results. Thus let 
$$\hat \cI_{\dm}: \cE_{+}^{m,\a} \to \bR^{+} \cup \{0\}$$
be a continuous approximation to $\frac{1}{\cO_{\dm}}$ with $\hat \cI_{\dm} \geq \frac{1}{\cO_{\dm}}$. 

  Moreover, via smooth local partitions of unity, one may convolve $\hat \cI_{\dm}$ with a 
smooth local mollifier to obtain a $C^{\infty}$ smooth approximation $\w \cI_{\dm}$ to $\hat \cI_{\dm}$ 
which satisfies  
\be \label{half}
\w \cI_{\dm} \geq \frac{1}{2} \cI_{\dm}.
\ee
We choose the approximation and smoothing above so that 
\be \label{zeroset}
\w \cI_{\dm} = 0 \ \ {\rm on} \ \ \w \cV^{m,\a}, 
\ee
where $\w \cV^{m,\a}$ is a proper open subset of  $\cV^{m,\a}$ containing the exterior Schwarzschild metrics 
\eqref{Sch}. The Schwarzschild metrics have $\cO_{\dm} = \infty$ so $\cI_{\dm} = 0$. 

  Next the function $H \to \min H$ is a continuous function on $C_{+}^{m-1,\a}(\dm)$, with values in $\bR^{+}$. 
The ratio
$$\frac{H}{\min H}: C_{+}^{m-1,\a}(\dm) \to \bR^{+},$$
is thus also continuous. Analogous to the above, let $min H$ be a smooth approximation to  $\min H$. 
The ratio
$$\frac{H}{min H}: C_{+}^{m-1,\a}(\dm) \to \bR^{+},$$  
is thus smooth. 

  Given this background, consider now the map
\be \label{wPi}
\w \Pi: \cE^{m,\a} \to Met^{m,\a}(S^{2})\times C_{+}^{m-1,\a}(S^{2}),
\ee
$$\w \Pi(g, u) = (\g, \w H),$$
where 
\be \label{wH}
\w H = H(1 + \frac{1}{min H}\w \cI_{\dm}) = H + \frac{H}{min H}\w \cI_{\dm}.
\ee
Thus we are altering the boundary map $\Pi_{B}$ by a scaling factor $\l = (1 + \frac{1}{min H}\w \cI_{\dm})$ 
on $H$. This choice of scaling factor is not unique or canonical any sense; it is a useful or convenient choice. 
While $(\g, H)$ depend only on the boundary data, the term $\w \cI_{\dm}$ approximating $\frac{1}{\cO_{\dm}}$ 
depends on the global geometry of $(M, g, u) \in \cE^{m,\a}$. 

  Observe that 
\be \label{piwpi}
\w \Pi = \Pi_{B}
\ee
in domains of $\cE_{+}^{m,\a}$ where $\w \cI_{\dm} = 0$, i.e.~on the domain $\w \cV^{m,\a}$ defined in 
\eqref{zeroset}. As noted above, the standard exterior $\{r \geq R\}$ regions of the Schwarzschild metric are in 
$\w \cV^{m,\a}$. Note also that the family of maps 
$$\w \Pi_{\d}(g, u) = (\g, \w H_{\d}),$$
where 
$$\w H_{\d} = H(1 + \frac{\d}{min H}\w \cI_{\dm}),$$ 
interpolates continuously between $\Pi_{B}$ and $\w \Pi$ as $\d$ ranges over the interval $[0,1]$. 

\begin{proposition} 
The map $\w \Pi$ is a smooth Fredholm map of Banach spaces, of Fredholm index $0$.
\end{proposition}

{\bf Proof:} The map $\w \Pi$ is a smooth map between Banach manifolds, by construction. To prove it is 
Fredholm, recall from the Introduction that the boundary map $\Pi_{B}$ is smooth and Fredholm, 
(as proved in \cite{An1}, \cite{AK}). The derivatives of $\w \Pi$ and $\Pi_{B}$ are closely related. One has 
\be \label{Dpi}
D\w \Pi(h, u') = (h^{T}, H'_{h}(1 + \frac{1}{min H}\w \cI_{\dm})) + (0, H(\frac{1}{min H}\w \cI_{\dm})').
\ee
Here $h^{T}$ is the restriction of the metric variation $h$ of $g$ to $\dm$, whle $H'_{h}$ is the induced 
variation of the mean curvature. The first term in \eqref{Dpi}, denoted $D_{1}$, differs from $D\Pi_{B}$ 
by a single fixed scale factor $(1 + \frac{1}{min H}\w \cI_{\dm})$ and since $D\Pi_{B}$ is Fredholm, 
it is clear that $D_{1}$ is also a Fredholm map of Banach spaces. For the second term in \eqref{Dpi}, denoted 
$D_{2}$, the derivative (at any given point $(M, g, u) \in \cE_{+}^{m,\a}$) of the smooth function 
$(\frac{1}{min H} \w \cI_{\dm})'$ is a bounded linear functional of $(h, u')$ and so $H(\frac{1}{min H} \w \cI_{\dm})'$ 
is a compact operator (in fact an operator of rank at most $1$). Since the sum of a Fredholm operator and a 
compact operator is Fredholm, this proves that $\w \Pi$ is Fredholm. 

  Exactly the same argument shows that the interpolating maps $\w \Pi_{\d}$ are also smooth and Fredholm. 
It is well-known that the Fredholm index (cf.~\eqref{ind0}) is deformation invariant, so that the Fredholm 
index of $D \w \Pi$ is independent of the point $(M, g, u) \in \cE^{m,\a}$ at which the linearization
is taken, in any component of $\cE_{+}^{m,\a}$. Similarly, the Fredholm index of $\w D\Pi_{\d}$ 
(at any base point) is independent of $\d$. Since the Fredholm index of $D\Pi_{B}$ is zero 
(everywhere) it follows that the index of $D\w \Pi$ is also zero everywhere. 

{\endproof} 

   The key to obtaining global information from a Fredholm map is to prove that the map is proper. 
Recall that a continuous map $F: X \to Y$ is called {\it proper} if, for any compact set $K \subset Y$, 
the inverse image $F^{-1}(Y)$ is compact in $X$. In the context of the map $\w \Pi$, this amounts to 
showing that control of the boundary data $(\g, \w H)$ (i.e.~data in a compact set of the target space) 
implies control of corresponding solutions $(M, g, u)$ realizing such boundary data. 

  The examples of boundary degeneration discussed following \eqref{PiB+} show that $\Pi_{B}$ is not proper. 

\begin{theorem}
The map $\w \Pi$ is proper. 
\end{theorem}

{\bf Proof:} Let $(M_{i}, g_{i}, u_{i})$ be a sequence in $\cE_{+}^{m,\a}$ such that $\w\Pi (M_{i}, g_{i}, u_{i})$ 
converges to a limit in $Met^{m,\a}(S^{2})\times C_{+}^{m-1,\a}(S^{2})$, so that 
\be \label{bcont}
\g_{i} \to \g_{\infty}, \ \ \w H_{i} \to \w H_{\infty},
\ee
in the $C^{m,\a}$ and $C^{m-1,\a}$ topologies respectively. We need to prove that a subsequence of 
$(M_{i}, g_{i}, u_{i})$ converges to a limit in $\cE_{+}^{m,\a}$. The basic idea is to use Proposition 2.1. 
For this we show that \eqref{bcont} implies that the hypotheses of Proposition 2.1 hold. (Of course $S$ 
is taken to be the set $\{(M, g_{i}, u_{i})\}$). 

  We first claim there is a constant $C_{0} < \infty$ such that (on the sequence $(M, g_{i}, u_{i})$) 
\be \label{obound}
\w \cI_{\dm} \leq C_{0} < \infty. 
\ee
Namely, since $H > 0$ in $\cE_{+}^{m,\a}$ one has 
$$\w H >  \frac{H}{\min H}\w \cI_{\dm} \geq \w \cI_{\dm},$$
since $\frac{H}{\min H} \geq 1$. Thus \eqref{bcont} implies \eqref{obound}.  
Via \eqref{half}, the bound \eqref{obound} implies that $\cO_{\dm} \geq c_{0} > 0$, for some constant $c_{0} > 0$.
  
  Next we show that \eqref{bcont} implies 
\be \label{Hconv}
H_{i} \to H_{\infty} \ \ {\rm and} \ \ H_{\infty} > 0,
\ee
so that $H_{i} \to H_{\infty} \in C_{+}^{m-1,\a}(S^{2})$. 
  
  To do this, we first prove that $\min H_{i}$ is bounded away from $0$; there exists $h_{0} > 0$ such that 
\be \label{hlow}
H_{i} \geq h_{0} > 0.
\ee
For if not, then 
\be \label{not}
\min H_{i} \to 0
\ee
(in a subsequence). Choose points $x_{i}$ such that $H(x_{i}) = \min_{x \in \dm} H_{i}(x)$. 
Then $\w H(x_{i}) = H(x_{i}) + \w \cI_{\dm}$ and since $\w H(x_{i})$ is bounded away from $0$ by \eqref{bcont}, 
it follows that there exists a constant $c_{0} > 0$ such that 
\be \label{obound2}
c_{0} \leq \w \cI_{\dm} \leq C_{0}. 
\ee
It then follows from \eqref{obound2} and \eqref{bcont} that the function 
$$\frac{H}{\min H} + H$$ 
is uniformly bounded away from $0$ and $\infty$ on the sequence $(M, g_{i}, u_{i})$. Since by assumption 
$\min H_{i} \to 0$, this implies that 
$$H_{i} \to 0$$
uniformly on $\dm$. (The geometry of $\dm$ is uniformly controlled also by \eqref{bcont}). By \eqref{obound} 
(i.e.~$\cO_{\dm} \geq c_{0} > 0$) and Remark 4.4 of \cite{AK}, a subsequence of $(M, g_{i}, u_{i})$ then converges 
in $\cE^{m,\a}$ to a limit $(M, g, u) \in \cE^{m,\a}$ with boundary data $(\g, H) = (\g_{\infty}, 0)$. (This is the extension 
of Proposition 2.1 to the closure $\bar \cE_{+}^{m,\a}$). A result of Miao, \cite{M}, which is a variant of the black hole uniqueness 
theorem, implies that $(M, g, u)$ is the Schwarzschild metric \eqref{Sch} for some mass $m$. In particular, if 
$\g_{\infty}$ in \eqref{bcont} is not a round metric of some fixed radius, then one already has a contradiction to 
\eqref{not}. The Schwarzschild metric has $\cO_{\dm} = \infty$ and is in the space $\w \cV^{m,\a}$ where 
$\w \cI_{\dm} = 0$. Since $\w \cV^{m,\a}$ is open, $(M, g_{i}, u_{i}) \in \w \cV^{m,\a}$ for $i$ sufficiently large. This 
contradicts \eqref{obound2} which then establishes \eqref{hlow}. 
 
  Thus \eqref{bcont} implies that $\min H_{i}$ is uniformly bounded away from $0$ and 
$\infty$ and $\w \cI_{\dm}$ is uniformly bounded above on the sequence $(M, g_{i}, u_{i})$. From this, 
it follows trivially from \eqref{bcont} that $H_{i} \to H > 0$ in $C_{+}^{m-1,\a}(S^{2})$ (in a subsequence). 
This proves \eqref{Hconv}. 

   The result then follows from Proposition 2.1.

{\endproof}

   A proper Fredholm map $F$ of index 0 between connected and separable Banach manifolds has a 
$\bZ_{2}$-valued degree, the Smale degree $deg_{\bZ_{2}}F$, given by 
\be \label{deg}
deg_{\bZ_{2}}F = \# \{F^{-1}(y)\} \ \ {\rm mod} \ 2,
\ee
where $y$ is a regular value of $F$, i.e.~for any $x \in F^{-1}(y)$ the linearization $(DF)_{x}$ is 
surjective. Here $\#$ denotes the cardinality; the set $\{F^{-1}(y)\}$ is a finite set of points, by the 
inverse function theorem and properness assumption.  Of course the degree is independent of the 
choice of regular value $y$. We refer to \cite{Sm} (or also \cite{AK} here) for further detail. 

  Apriori, the separable Banach manifold $\cE_{+}^{m,\a}$ may not be connected. Let then $\cE_{+,0}^{m,\a}$ 
be the component of $\cE_{+}^{m,\a}$ containing the exterior Schwarzschild metrics as described following \eqref{Sch}. 
This is the same as the component containing the flat static vacuum solutions $(M, g_{Eucl}, 1)$ where $M$ the 
the region exterior to a smooth compact convex subset $B \subset \bR^{3}$. 

   The map $\w \Pi$ restricts to a map 
\be \label{conn}
\w \Pi: \cE_{+,0}^{m,\a} \to Met^{m,\a}(S^{2})\times C_{+}^{m-1,\a}(S^{2}).
\ee
By Proposition 2.3 and Theorem 2.4, $\w \Pi$ is a proper Fredholm map of index 0 between connected and 
separable Banach manifolds. Hence 
$$deg_{\bZ_{2}} \w \Pi \in \bZ_{2},$$
is well-defined. 

\begin{proposition}
For the map $\w \Pi$ in \eqref{conn}, one has  
\be \label{deg1}
deg_{\bZ_{2}}\w \Pi = 1.
\ee
In particular $\w \Pi$ is surjective onto $Met^{m,\a}(S^{2})\times C_{+}^{m,\a}(S^{2})$. 
\end{proposition}

{\bf Proof:}  The proof is a simple modification of the proof of Theorem 5.1 of \cite{AK}. On the curve of 
Schwarzschild metrics \eqref{Sch} with $\dm$ taken at a locus with $r$ fixed, one has 
$$(\g, H) = (\g_{S^{2}(r)}, \frac{2}{r}\sqrt{1 - \frac{2m}{r}}).$$
Thus any constant curvature metric $\g_{S^{2}(r)}$, and any constant $H_{0} \in (0,\infty)$ are realized in 
$Im \Pi_{B}$. Since $\cI_{\dm} = 0$ on the Schwarzschild curve (so that \eqref{piwpi} holds) the same 
boundary data satisfy $(\g_{S^{2}(r)}, H_{0}) \in Im(\w \Pi)$. 

 Now fix $r = 1$ and conversely consider boundary data as above $(\g_{S^{2}(1)}, t_{i})$ for $\w \Pi$, 
with $\w H_{t_{i}} = t_{i} \to 0$ as $t_{i} \to 0$. Let $(M, g_{i}, u_{i})$ be any sequence of points in 
$(\w \Pi)^{-1}(\g_{S^{2}(1)}, t_{i})$. On such a sequence one has  
$$H(1 + \frac{1}{min H}\w \cI_{\dm}) = t_{i} \to 0,$$
which implies that 
$$H = const \to 0 \ \ {\rm and} \ \ \cI_{\dm} \to 0.$$
It then follows then as above in the proof of Theorem 2.4, namely from [3, Remark 4.4] (see also the beginning 
of the proof of [3, Theorem 5.1]) that $(M, g_{i}, u_{i})$ converges again to the Schwarzschild metric 
of mass $m = \frac{1}{2}$. 

   Now suppose 
\be \label{deg0}
deg_{\bZ_{2}}\w \Pi = 0.
\ee
By the Sard-Smale theorem \cite{Sm}, regular values of $\w \Pi$ are open and dense. Choose a sequence of regular 
values $(\g_{i}, H_{i})$ in $Met^{m,\a}(S^{2})\times C_{+}^{m-1,\a}(S^{2})$ converging to $(\g_{S^{2}(1)}, 0)$. 
For any $i$, the finite set $(\w \Pi)^{-1}(\g_{i}, H_{i})$ (if non-empty) consists of at least two distinct static 
vacuum solutions $(M, g_{i}^{1}, u_{i}^{1})$, $(M, g_{i}^{2}, u_{i}^{2})$. As above, it follows that both sequences 
converge to the Schwarzschild metric \eqref{Sch} with $m = \frac{1}{2}$. 

  Thus, as discussed following \eqref{zeroset} (and in the proof of Theorem 2.4 above), 
$$\w \cI_{\dm}(g_{i}^{j}) = 0,$$ 
for $j = 1,2$ and $i$ sufficiently large, so that $\w \Pi = \Pi_{B}$ at and near each $g_{i}^{j}$. The rest of the 
proof of [3, Theorem 5.1] now applies verbatim to give a contradiction to \eqref{deg0}, which establishes 
\eqref{deg1}. 

   It is standard that any map of degree one is necessarily surjective. (Any point not in the image of $F$ is a
regular value of $F$). 
 
{\endproof}

  Proposition 2.5 leads trivially to the proof of Theorem 1.2. 

\medskip 
  
  {\bf Proof of Theorem 1.2.}
  
   Let $\pi: Met^{m,\a}(S^{2})\times C_{+}^{m-1,\a}(S^{2}) \to Met^{m,\a}(S^{2})\times \cC_{+}^{m-1,\a}(S^{2})$
be the projection to the orbit space of the $\bR^{+}$ action as following \eqref{equiv}. The map $\pi$ is clearly 
surjective, and hence by Proposition 2.5, $\pi \circ \w \Pi: \cE_{+}^{m,\a} \to 
Met^{m,\a}(S^{2})\times \cC_{+}^{m-1,\a}(S^{2})$ is also surjective. Since 
\be \label{Pipi}
\Pi = \pi \circ \w \Pi,
\ee
this proves the result. 
  
{\endproof}

{\bf Proof of Theorem 1.1.}

   Let $(\g, H) \in Met^{\infty}(S^{2})\times C_{+}^{\infty}(S^{2})$ be smooth boundary data. Then of course $(\g, H)$ is 
in the little H\"older space $Met^{m,\a}(S^{2})\times C_{+}^{m,\a}(S^{2})$, for any $m$, $\a$. By Theorem 1.2, there 
exists a static vacuum solution $(M, g, u) \in \cE_{+,0}^{m,\a} \subset \cE_{+}^{m,\a}$ realizing the boundary data 
$(\g, [H])$, so the mean curvature of $(M, g, u)$ at $\dm$ is given by $\l H$ for some $\l > 0$. The data $(g, u)$ are 
$C^{m,\a}$ up to $\dm$. Recall from the Introduction that the static vacuum Einstein equations with boundary data 
$(\g, H)$ form an elliptic boundary value problem (in a suitable gauge). By elliptic boundary regularity, it follows 
that $(g, u)$ are $C^{\infty}$ up to $\dm$, so that $(M, g, u) \in \cE_{+,0}^{\infty}$. This completes the proof. 

{\endproof}

\section{Discussion}
\setcounter{equation}{0} 
  
  Consider the map $\Pi$ in \eqref{Pi}. By the Sard-Smale theorem \cite{Sm}, regular values 
are open and dense in the target space $Met^{m,\a}(S^{2})\times \cC_{+}^{m-1,\a}(S^{2})$; as in \eqref{deg}, 
$y$ is a regular value of $F$ if $(DF)_{x}$ is surjective, for each $x \in F^{-1}(y)$. 

   Let $(\g, [H])$ be any regular value of $\Pi$. Since $\w \Pi$ is proper and $\pi$ as in \eqref{Pipi} 
has 1-dimensional fibers, it follows that 
\be \label{curves}
\Pi^{-1}(\g, [H]),
\ee
is a finite collection of smooth curves in $\cE_{+}^{m,\a}$. This corresponds to \eqref{ind1}, that the 
Fredholm index of $\Pi$ in \eqref{Pi} is one. Let $(M_{t}, g_{t}, u_{t})$ denote 
any one such curve. The induced boundary data are $(\g, H_{t})$ with $H_{t} = \l_{t}H$, for some 
fixed function $H$, so that $\l_{t} \in \bR^{+}$ depends smoothly on $t$.   

\medskip 
  
  {\bf Question:} Is $\l_{t}$ surjective onto $\bR^{+}$?

\smallskip 
  
     This is the question of whether $\Pi_{B}$ in \eqref{PiB} is surjective onto the $\bR^{+}$ orbits of 
$H$; it is surjective onto the orbit space by Theorem 1.2. If this is the case (which would be surprising in 
the author's opinion) then of course the Bartnik boundary map $\Pi_{B}$ is essentially surjective; it has at 
least open and dense range. However, as discussed in the Introduction, one does not expect $\Pi_{B}$ to be 
injective or have a well-behaved inverse. Consider for example the curve of metrics $(M_{t}, (g_{Eucl})_{t}, 1)$ 
in $\bR^{3}$ where $\dm_{t}$ converges to an immersed sphere of positive mean curvature in $\bR^{3}$ as 
$t \to 0$. The limit configuration $(M_{0}, (g_{Eucl})_{0}, 1)$ makes sense, and has a well-defined boundary 
value $(\g_{0}, H_{0}) \in Met^{m,\a}(S^{2})\times C_{+}^{m-1,\a}(S^{2})$, but $(M_{0}, (g_{Eucl})_{0}, 1)$ 
is not a solution in the space $\cE_{+}^{m,\a}$. On the other hand, Theorem 1.2 implies there is a $\l > 0$ 
such that $(\g_{0}, \l H_{0})$ is realized as boundary data of a solution $(M, g, u) \in \cE_{+}^{m,\a}$. 
  
\medskip 

 Next, an immediate consequence of Theorem 1.2 is that any $C^{m,\a +}$ metric $\g$ on $S^{2}$ may 
be extended to a smooth asymptotically flat metric of zero scalar curvature on $M$, with mean 
curvature $H > 0$ on $\dm$ prescribed up to a multiplicative constant. Modulo the scalar 
indeterminacy on $H$, this extends previous work of Shi-Tam \cite{ST} who showed that 
boundary data $(\g, H)$ can be extended to scalar-flat asymptotically flat metrics provided 
the Gauss curvature $K_{\g} > 0$. The work of Shi-Tam is based on the method of Bartnik 
\cite{B3} on constructing asymptotically flat metrics of prescribed (non-negative) scalar 
curvature with given boundary data $(\g, H)$. 

\medskip 

   Finally we discuss issues related to the mass of solutions in $\cE_{+}^{m,\a}$. Let $m$ denote the 
ADM mass of a static vacuum solution $(M, g, u)$. As is well-known, $m$ is given as a Komar-type 
integral (mass), 
$$m = \frac{1}{4\pi}\int_{\dm}N(u)dv_{\g}.$$
In particular, the mass 
$$m: \cE_{+}^{m,\a} \to \bR$$
is a smooth function on $\cE_{+}^{m,\a}$. Clearly the space $\cE_{+}^{m,\a}$ contains solutions with both
positive and negative mass (e.g.~the Schwarzschild curve). Thus $\cH_{0} = \{m = 0\}$ is a hypersurface in 
$\cE_{+}^{m,\a}$; it seems likely that $\cH_{0}$ is a smooth hypersurface, but this remains to be proved. 

  Let $(M_{t}, g_{t}, u_{t})$ be a smooth curve component of $\Pi^{-1}(\g, [H])$ as in \eqref{curves}, so that  
$m_{t}$ is a smooth function of $t$. 

\medskip 
   
   {\bf Question:} Is $\l_{t}$ a function of $m_{t}$ or vice versa? Is $\l_{t}$ monotone in $m_{t}$? 

\smallskip 
   
    Note this is the case for the curve of Schwarzschild metrics $g_{m}$ described following \eqref{Sch}; 
here $\l_{t}$ takes all values in $(0, \infty)$ and is monotone decreasing as a function of $m_{t} \in 
(-\infty, \frac{1}{2})$. 
  
  \medskip 
  
  {\bf Question:} Are there criteria on $(\g, H)$ which guarantee that $m > 0$ or $m \geq 0$? 
  
\smallskip 
   
    This question is motivated by the original context of the Bartnik conjecture and the definition of 
the Bartnik quasi-local mass \cite{B1}, \cite{B2}, which we briefly recall. Let $(\hat M, \hat g)$ be a 
complete asymptotically flat manifold of non-negative scalar curvature and $\Omega$ a bounded 
smooth domain in $\hat M$ with $H > 0$ at $\partial \Omega$. Let $M = \hat M \setminus \Omega$, 
assumed to be diffeomorphic to $\bR^{3}\setminus B$. Then $\dm = \partial \Omega \simeq S^{2}$ 
with induced Bartnik boundary data $(\g, H)$. A special case of the Bartnik conjecture is that $(\g, H)$ are 
boundary data of a unique solution $(M, g, u) \in \cE_{+}^{m,\a}$. By the results of \cite{ST} and the 
positive mass theorem, $m \geq 0$ on $(M, g, u)$. In this context, one would like to characterize or understand 
the space of solutions in $\cE_{+}^{m,\a}$ which have $m > 0$ or $m \geq 0$ in terms of the 
boundary data $(\g, H)$.

   A result of Jauregui \cite{J} shows that, given such data $(\g, H)$ on $S^{2}$, and assuming 
$K_{\g} > 0$, there is a maximal value $\l_{0} < \infty$ (depending on $(\g, H)$) such that 
$(\g, \l H)$ has an in-filling 3-ball $B^{3}$ with a metric of non-negative scalar curvature for all $\l < \l_{0}$, 
but no such in-filling exists for $\l > \l_{0}$. See also \cite{JMT} for a refinement of this result. 

  It would be interesting to know if such a criterion holds for the mass, e.g.~whether $m \geq 0$ if any only if 
$\l \leq \l_{0}$, for some $\l_{0}$. This is true on the Schwarzschild curve, where $m \geq 0$ if and only 
if $H \leq 2$ for $\g = \g_{S^{2}(1)}$.

\bibliographystyle{plain}

\begin{thebibliography}{WWW}

\footnotesize


\bibitem [1]{An1} M.~T.~Anderson, {\it On boundary value problems for Einstein metrics}, 
Geom. \& Topology, {\bf 12}, (2008), 2009-2045, arXiv:math/0612647

\bibitem[2]{An2} M.~T.~Anderson, {\it Static vacuum Einstein metrics on bounded domains}, 
Annales Henri Poincar\'e, {\bf 16}, (2015), 2265-2302, arXiv:1305.1540.

\bibitem[3]{AK} M.~T.~Anderson and M. Khuri, {\it On the Bartnik extension problem for static vacuum 
Einstein metrics}, Classical \& Quantum Gravity, {\bf 30}, (2013), 125005, arXiv:0909.4550.  


\bibitem[4]{B1} R. Bartnik, {\it New definition of quasi-local mass}, Phys. Rev. Lett., {\bf 62},
(1989), 2346-2348.

\bibitem[5]{B2} R. Bartnik, {\it Energy in general relativity}, Tsing Hua lectures on geometry
and analysis, (Hsinchu, 1990-91), 5-27, International Press, Cambridge, MA, 1997.

\bibitem[6]{B3} R. Bartnik, {\it Quasi-spherical metrics and prescribed scalar curvature}, Jour. 
Differential Geom., {\bf 37}, (1993), 31-71. 

\bibitem[7]{B4} R. Bartnik, {\it Mass and 3-metrics of non-negative scalar curvature},
Proc. Int. Cong. Math., vol II, Beijing (2002), 231-240, Higher Ed. Press, Beijing, 2002, 
arXiv:math/0304259.

\bibitem[8]{BFT} R. Bonic, J. Frampton and A. Tromba, {\it $\Lambda$-manifolds}, 
Jour. Functional Analysis, {\bf 3}, (1969), 310-320. 

\bibitem[9]{D} P. Delano\"e, {\it Local inversion of elliptic problems on compact manifolds}, 
Math. Japonica, {\bf 35}, (1990), 679-692. 

\bibitem[10]{FM} R. Fry and S. McManus, {\it Smooth bump functions and the geometry of Banach 
spaces: A brief survey}, Expositiones Math., {\bf 20}, (2002), 143-183. 

\bibitem[11]{J} J.~L.~Jauregui, {\it Fill-ins of non-negative scalar curvature, static metrics, and quasi-local mass}, 
Pacific Jour. Math., {\bf 261}, (2013), 417-444, arXiv:1106.4339. 

\bibitem[12]{JMT} J.~L.~Jauregui, P. Miao and L-F.~Tam, {\it Extensions and fill-ins with non-negative scalar 
curvature}, Classical \& Quantum Gravity, {\bf 30}, (2013), 195007, arXiv:1304.0721

\bibitem[13]{M} P. Miao, {\it A remark on boundary effects in static vacuum initial data sets}, 
Classical \& Quantum Gravity, {\bf 22}, (2005), L53-59. 

\bibitem[14]{Mu} J. M. Munkres, {\it Topology}, $2^{\rm nd}$ Edition, Prentice-Hall, (2000). 

\bibitem[15]{O} D. Opela, {\it Spaces of functions with bounded and vanishing mean oscillation}, 
in: {\it Function Spaces, Differential Operators and Nonlinear Analysis}, (Teistungen, 2001), Birkh\"auser Verlag, 
Basel, (2003), 403-413. 


\bibitem[16]{ST} Y. Shi and L-F.~Tam, {\it Positive mass theorem and the boundary behaviors of compact 
manifolds with non-negative scalar curvature}, Jour. Differential Geom., {\bf 62}, (2002), 79-125, 
arXiv:math/0301047. 


\bibitem[17]{Sm} S. Smale, {\it An infinite dimensional version of Sard's theorem},
Amer. Jour. Math., {\bf 87}, (1965), 861-866.



\end{thebibliography}

\end{document}